\documentclass[12pt]{article}

\usepackage[affil-it]{authblk}
      \usepackage{latexsym}
         \usepackage[reqno, namelimits, sumlimits]{amsmath}
         \usepackage{amssymb, amsfonts}
         \usepackage{amsthm}

 \newtheorem{theorem}{Theorem}[section]
 
 \newtheorem{definition}[theorem]{Definition}

 \newtheorem{pro}[theorem]{Proposition}

\title{A note on weak solutions to the Navier-Stokes equations that are locally  
in $L_\infty(L^{3,\infty})$}

\author{ G Seregin
  \thanks{ seregin@maths.ox.ac.uk; }
  }
\affil{OxPDE, Mathematical Institute, University of Oxford, Oxford,UK, and 
St Petersburg Department of V A Steklov Mathematical Institute, Russia}
\date{ \today}

\begin{document}
\maketitle

\centerline{Dedicated to Nina Nikolaevna Uraltseva }

\begin{abstract}
The aim of the note is to proof a regularity result for weak solutions to the Navier-Stokes equations that are locally  
in $L_\infty(L^{3,\infty})$. It reads that, in a sense, the number of singular points at each time is at most finite. Our note is inspired by the paper of H. J. Choe, J. Wolf, M. Yang \cite{CWY}.

\end{abstract}

\setcounter{equation}{0}
\section{Inroduction}

Our note is very much motivated by the paper \cite{CWY}. The authors of  \cite{CWY} consider  a weak solution to the Cauchy problem for the Navier-Stokes equations with $L_2$-initial data  under the additional assumoption that it is bounded  in time with values in the weak Lebesgue space $L^{3,\infty}(\mathbb R^3)$. They show that, at each instance of time, there exists  at most a finite number of singular points.

What we would like  is to extend this result to the local setting, including considerations near a flat part of the  boundary,  and to the standard notion of suitable solutions, see \cite{CKN}, \cite{Lin}, and \cite{LadSer1999}. Our proof seems to be shorter and straightforword.

Interior and bounded regularity will be analyised separately. 
Let us start with the interrior case.

Consider a suitable weak solution $v$ and $q$ in $Q_T=\Omega\times ]0,T[$, where $\Omega$ is a domain $\mathbb R^3$. The corresponding definition is due to F-H Lin, see \cite{Lin} and Definition \ref{sws} of the this paper.
 It differs slightly from the original one, introduced by Caffarelli-Kohn-Nirenberg in \cite{CKN}, just by a more convienient class for the pressure field.
\begin{definition}\label{sws} We say that a pair $v$ and $q$ is a suitable weak solution to the Navier-Stokes equations in $Q_T$ if:
\begin{equation}
	\label{class}
	v\in L_{2,\infty}(Q_T),\qquad v\in L_{2,\infty}(Q_T),\qquad q\in L_\frac 32(Q_T);
	\end{equation}

the pair $v$ and $q$ satisfies the Navier-Stokes equations	
	\begin{equation}
		\label{NSE}
		\partial_tv+v\cdot\nabla v-\Delta v=-\nabla q,\qquad {\rm div}\,v=0
	\end{equation}
in $Q_T$ in the sense of distributions;	

for a.a. $t\in ]0,T[$, the local energy inequality
\begin{equation}
	\label{locenineq}
\int\limits_{\Omega}|v(x,t)|^2\varphi(x,t)dx+2\int\limits^t_0\int\limits_\Omega |\nabla v|^2\varphi dx dt'\leq \end{equation}
$$\leq \int\limits^t_0\int\limits_\Omega(|v|^2(\partial_t\varphi+\Delta \varphi)+v\cdot\nabla \varphi(|v|^2+2q))dx dt'
$$	
holds for all non-negative test functions $\varphi \in C^\infty_0(\Omega \times ]0,2T[)$.
\end{definition}

Our basic additional assumpton is that 
\begin{equation}
	\label{addass}
	\|v\|_{ L_\infty(0,T;L^{3,\infty}(\Omega))}\leq M<\infty.
\end{equation}
In fact, it implies the following: one can select a representative of the function $t\to v(\cdot,t)$ so that 
\begin{equation}
	\label{conaddass}
	\sup\limits_{0<t\leq T}\|v(\cdot,t)\|_{L^{3,\infty}(\Omega)}\leq M.\end{equation}
Indeed, fix a representative for $v$ such that the set of all singular points has zero 1D parabolic Hausdorff measure. Hence, for each time $0<t_0\leq T$, the set of singular points $(x,t_0)$ has zero 1D Hausdorff measure. As it has been shown in \cite{LadSer1999}, the function $z=(x,t)\to v(z)$ is H\"older continuous in a parabolic vicinity of each regular point $(x,t_0)$. So, the following is true:
$$v(x,t)\to v(x,t_0)$$
for a.a. $x\in \Omega$ as $t\to t_0$ and $t<t_0$. Then, selecting a sequence of times $t_k<t_0$ such that 
 $$\|v(\cdot,t_k)\|_{L^{3,\infty}(\Omega)}:=\sup\limits_{\alpha>0}\alpha|\{x\in \Omega:\,\,|v(x,t_k)|>\alpha\}|^\frac 13\leq M,
$$
observe that
$$\liminf\limits_{k\to \infty}\|v(\cdot,t_k)\|_{L^{3,\infty}(\Omega)}\geq\|v(\cdot,t_0)\|_{L^{3,\infty}(\Omega)}.
$$

It is important to notice that condition (\ref{addass}) provides the existence of non-trivial limit solutions that are arising from rescaling procedure around a singular point.

A local version of the main result of the paper 
\cite{CWY} can be proved with the help of an idea  from the paper \cite{Seregin1999}.
\begin{theorem}
\label{final} Let $v$ and $q$ be a suitable weak solution to the Navier-Stokes equations in $Q_T$. Assume that $v$ satisies condition (\ref{addass}). Then,
for any subdomain $\Omega_1\Subset\Omega$,  there exists  at most a finite number of singular points in the set $\{(x,T):\,\,x \in \overline\Omega_1\}$.
\end{theorem}



To prove Theorem \ref{final}, we need intermediate statements that might be interesting themselves. In order to describe them, let us introduce the following scale invariant quantities:

$$A(v,r;z_0):=\frac 1r\sup\limits_{t_0-r^2<t<t_0}\int\limits_{B(x_0,r)}|v(x,t)dx,\qquad C(v,r;z_0)=\frac 1{r^2}\int\limits_{Q(z_0,r)}|v|^3dz,
$$
$$
E(v,r;z_0)=\frac 1{r}\int\limits_{Q(z_0,r)}|\nabla v|^2dz,\qquad K(v,r;z_0)=\frac 1{r}\int\limits_{Q(z_0,r)}|v|^4dz,$$
$$ D(r)=\frac 1{r^2}\int\limits_{Q(r)}|p|^\frac 32dz, \qquad
D_0(q,r;z_0)=\frac 1{r^2}\int\limits_{Q(z_0,r)}|q-[q]_{B(x_0,r)}|^\frac 32dz,
$$
where $z_0=(x_0,t_0)$, $Q(z_0,r)=B(x_0,r)\times ]t_0-r^2,t_0[$, $B(x_0,r)=\{|x-x_0|<r\}$, $[q]_{B(x_0,r)}(t)$ is the mean value of the function $x\to q(x,t)$ over the ball $B(x_0,r)$. 
Also, let us abbreviate: $B(r)=B(0,r)$, $B=B(1)$, $Q(r)=Q(0,r)$, $Q=Q(1)$,
$A(v,r)=A(v,r;0)$, etc.

The following proposition is a local version of the main regularity result of the paper \cite{CWY}.

\begin{pro}\label{mainlemma} 
Let $v$ and $q$ be a suitable weak solution to the Navier-Stokes equations in $Q$, satisfying additional assumptions:
\begin{equation}
	\label{energy}
	D_0(q,1)+E(v,1)\leq N
	\end{equation}
and	
\begin{equation}
	\label{weakl3}
	\|v\|_{L_\infty(-1,0;L^{3,\infty}(B))}\leq M.
\end{equation}
There exists a positive number $\varepsilon
<\frac 14$, depending on $N$ and $M$ only, such that if, for some $0<r\leq \frac 12$, 
\begin{equation}
	\label{maincond}
	\frac 1{r^3}|\{x\in B(r): |v(x,0)|>\frac \varepsilon r\}|\leq\varepsilon,
\end{equation}
then 
\begin{equation}
	\label{boundedness}
	v\in L_\infty(Q(\varepsilon r)).\end{equation}
\end{pro}
In our further considerations,  a scaled version of Proposition \ref{mainlemma} is going to used. 
\begin{pro}\label{mainlemmasc} 
Let $v$ and $q$ be a suitable weak solution to the Navier-Stokes equations in $Q(z_0,R)$, satisfying additional assumptions:
\begin{equation}
	\label{energysc}
	D_0(q,R;z_0)+E(v,R:z_0)\leq N
\end{equation}
and	
\begin{equation}
	\label{weakl3sc}
	\|v\|_{L_\infty(t_0-R^2,t_0;L^{3,\infty}(B(x_0,R)))}\leq M.
\end{equation}
If, for some $0<r\leq \frac 12R$, 
inequality 
\begin{equation}
	\label{maincondscale}
	\frac 1{r^3}|\{x\in B(x_0,r): |v(x,t_0)|>\frac \varepsilon r\}|\leq\varepsilon
\end{equation}
holds, then $v\in L_\infty(Q(z_0,\varepsilon r))$.
\end{pro}

Next, let us discuss local regularity up to a flat part of the boundary. To formulate the corresponding results,  the specific notation is needed:
$$B^+(x_0,r)=B(x_0,r)\cap \{x_3>x_{03}\},\qquad Q^+(z_0,r)=B^+(x_0,r)\times ]t_0-r^2,t_0[.$$
 For $x_0=0$, abbreviations 
$B^+(r)=B^+(0,r)$, $B^+=B^+(1)$ are exploited. 

Now, the definition of suitable week solutiuons to the problem
\begin{equation}
	\label{bnse}
	\partial_tv+v\cdot\nabla v -\Delta v=-\nabla q,\qquad{\rm div}\,v=0
\end{equation}
in $Q^+$ and 
\begin{equation}
	\label{dirichlet}
	v(x',t)=0
\end{equation}
for all $-1<t<0$ ans for all $|x'|<1$, where
$$x'=(x_1,x_2,0)$$
for $x=(x_1,x_2,x_3)$, is as follows, see \cite{Seregin}.
\begin{definition}
	\label{swsbc}
We say that a pair $v$ and $q$ is a suitable weak solution to the Navier-Stokes equations in $Q^+$ if:
\begin{equation}
	\label{classbc}
	v\in L_{2,\infty}(Q^+),\qquad v\in L_{2,\infty}(Q^+),\qquad q\in L_\frac 32(Q^+);
	\end{equation}

the pair $v$ and $q$ satisfies  (\ref{bnse}) in the sense of distributions and $v$ satifies boundary condition (\ref{dirichlet});

for a.a. $t\in ]-1,0[$, the local energy inequality
\begin{equation}
	\label{locenineqbc}
\int\limits_{B^+}|v(x,t)|^2\varphi(x,t)dx+2\int\limits^t_{-1}\int\limits_{B^+} |\nabla v|^2\varphi dx dt'\leq \end{equation}
$$\leq \int\limits^t_{-1}\int\limits_B^+(|v|^2(\partial_t\varphi+\Delta \varphi)+v\cdot\nabla \varphi(|v|^2+2q))dx dt'
$$	
holds for all non-negative test functions $\varphi \in C^\infty_0(B \times ]-1,1[)$.
\end{definition}
Here, our main assumtion remain the same:
\begin{equation}
	\label{mabc}
	\|v\|_{L_\infty(-1,0;L^{3,\infty}(B^+))}\leq M<\infty.
\end{equation}	Arguing as above, one can show that $\|v(\cdot,t)\|_{L^{3,\infty}(B^+)}\leq M$ for all $t\in ]-1,0]$.

A boundary version of  our main result reads the folowing.
\begin{theorem}
	\label{finalbc}
Let $v$ and $q$ be a suitable weak solution to the Navier-Stokes equations in $Q^+$. Assume that $v$ satisies condition (\ref{mabc}). Then,
for any $r_0\in ]-1,0[$,  there exists  at most a finite number of singular points in the set $\{(x,0):\,\,x \in \overline{B}^+(r_0)\}$.
\end{theorem}

In order state the auxiliary results, let us define similar  scale invariant quantities, for example,
$$A^+(v,r;z_0)=\sup\limits_{t_0-r^2<t<t_0}\frac 1r\int\limits_{B^+(x_0,r)}|v(x,t)|^2dx,\, E^+(v,r;z_0)=\frac 1r\int\limits_{Q^+(z_0,r)}|\nabla v|^2dz,
$$ and so on. In addition, we introduce two other  pressure quantities:
$$D_2^+(q,r;z_0)=\frac 1{r^\frac {13}8}\int\limits^{t_0}_{t_0-r^2}\Big(\int\limits_{B^+(x_0,r)}|\nabla q|^\frac {12}{11}dx\Big)^\frac {11}8dt$$
and 
$$D_2(q,r;z_0)=\frac 1{r^\frac {13}8}\int\limits^{t_0}_{t_0-r^2}\Big(\int\limits_{B(x_0,r)}|\nabla q|^\frac {12}{11}dx\Big)^\frac {11}8dt.$$

Without loss of generality, one may assume that the suitable weak solution in  Definition \ref{swsbc} satisfies the additional condition
$D_2^+(q,1)=D_2^+(q,1;0)<\infty$.

Now,  an analog of Proposition \ref{mainlemma} can be stated as follows. 
\begin{pro}
	\label{mainlemmabc}
	Let $v$ and $q$ be a suitable weak solution to the Navier-Stokes equations in $Q^+$ in the sense of Definition \ref{swsbc}. Assume that it satisfies assumption (\ref{mabc}) and 
	\begin{equation}
		\label{energybc}
		D^+_2(q,1)+E^+(v,1)\leq N<\infty.
	\end{equation}
	There exists a positive constant $\varepsilon<\frac 14$, depending only on $N$ and $M$ only, such that, if,  for  some $0<r\leq 1/2$,
	\begin{equation}
		\label{maincondscalebc}
		\frac 1{r^3}|\{x\in B^+(r):\,|v(x,0)|>\frac \varepsilon r\}|\leq\varepsilon,
	\end{equation}
	then $v\in L_\infty(Q^+(\varepsilon r))$.
\end{pro}
 The scale version of Proposition \ref{mainlemmabc} reads the following.
\begin{pro}
	\label{mainlemmabcsc}
	Let $v$ and $q$ be a suitable weak solution to the Navier-Stokes equations in $Q^+(R)$ in the sense of Definition \ref{swsbc}. Assume that it satisfies assumption (\ref{mabc}) and 
	\begin{equation}
		\label{energybcsc}
		D^+_2(q,R)+E^+(v,R)\leq N<\infty.
	\end{equation}
	There exists a positive constant $\varepsilon<1/4$, depending only on $N$ and $M$ only, such that, if,  for  some $0<r\leq 1/2R$, inequality (\ref{maincondscalebc})
holds,	then $v\in L_\infty(Q^+(\varepsilon r))$.
\end{pro}

\setcounter{equation}{0}
\section{Proof of Theorem \ref{final}}

Let us fix an arbitrary subdomain $\Omega_1\Subset\Omega$ and let $\delta ={\rm dist}(\Omega_1,\partial\Omega)>0$
and $2R_\star=\min(\delta/2,\sqrt T)$. 

It is easy to verify that two  inequalities
$$A(v,r;z_0)\leq c\|v\|^2_{L_\infty(t_0-r^2,t_0;L^{3,
\infty}(B(x_0,r))}\leq cM^2$$
and 
$$K(v,r;z_0)\leq cM^2(E(v,r;z_0)+A(v,r;z_0)) $$
hold provided $Q(z_0,r)\subset Q_T$. Having those inequalities in hands and estimates for the energy scale invariant quantities proved in \cite{Seregin2006}, see Lemma 1.8, and in \cite{LadSer1999}, see Lemma 5.3, one can state that, for all $z_0=(x_0,T)$ with $x_0\in \Omega_1$, the following is true:
$$
	\sup\limits_{0<r<R_\star}A(v,r;z_0)+\sup\limits
_{0<r<R_\star}C^\frac 34(v,r;z_0)+\sup\limits_{0<r<R_\star}E(v,r;z_0)+
$$
\begin{equation}
	\label{scaleinv}+\sup\limits_{0<r<R_\star}K(v,r;z_0)+\sup\limits_{0<r<R_\star}D(q,r;z_0)
	\end{equation}
	$$\leq c(M)(D(q,R_\star;z_0)+E(v,R_\star;z_0)+1)\leq$$
	$$\leq  c(M,R_\star,
	 \|\nabla v\|_{L_2(Q_T)},\|q\|_{L_\frac 32(Q_T)})=:N.
	$$
The number $\varepsilon(M,N)$ of Proposition \ref{mainlemma} can be determined as numbers $M$ and $N$ are known.

Let $S$ be a set of all singular points of $v$ in  
$\{(x_0,T):\,\, x_0\in \Omega_1\}$. Assume that it contains more than $M^3\varepsilon^{-4}$ elements.  Letting $P=[ M^3\varepsilon^{-4}]+1$, one can find    $P$ different singular points $(x_k,T)$,  $k=1,2,...,P$, of the set $S$. Then, pick up  a positive number $R<R_*$  such that $B(x_k,R)\cap B(x_l,R)=\emptyset$ if $k\neq l$, $k,l=1,2,...,P$. According to Proposition \ref{mainlemmasc}, 
for  all $r\in ]0, 1/2R]$, the following should be true:
$$\varepsilon\leq \frac 1{r^3}|\{x\in B(x_k,r):|v(x,T)|>\frac \varepsilon r\}|$$
for all $k=1,2,...,P$. Now, we let $r=r_0=1/2R$ and, after summation over $k$, we arrive at the following inequality
$$P\varepsilon\leq \sum\limits_{k=1}^P\frac 1{r_0^3}|\{x\in B(x_k,r_0):|v(x,T)|>\frac \varepsilon {r_0}\}|=
$$
$$=\frac 1{r_0^3}|\{x\in\bigcup\limits^P_{k=1}B(x_k,r_0):|v(x,T)|>\frac\varepsilon{r_0}\}|\leq 
$$$$
\leq \frac 1{r_0^3}|\{x\in \Omega:|v(x,T)|>\frac \varepsilon{r_0}\}|\leq \frac 1{\varepsilon^3}\|v(\cdot,T)\|^3_{L^{3,\infty}(\Omega)}\leq \frac {M^3}{\varepsilon^3}.
$$
The latter inequality implies that
$P\leq M^3\varepsilon^{-4}<P$. It is a contraduction. The theorem is proved.

\setcounter{equation}{0}
\section{Proof of Theorem \ref{finalbc}}

Let us first prove that the number of singular points of $v$ in the set $b(r_0)\times \{t=0\}$, where
$b(r_0)=\{x\in \mathbb R^3: \,x=x', |x'|\leq r_0\}$, is finite.

We let $2R_*=(1-r_0)/2$. Our further arguments  are very similar to ones used in the previous chapter, see \cite{Seregin}.  Indeed, for all space-time points $z_0=(x_0,0),$ where $x_0\in b(r_0)$, we have
$$
	\sup\limits_{0<r<R_\star}A^+(v,r;z_0)+\sup\limits_{0<r<R_\star}C^+(v,r;z_0)+\sup\limits_{0<r<R_\star}E^+(v,r;z_0)+
$$
\begin{equation}
	\label{scaleinvbc}+\sup\limits_{0<r<R_\star}K^+(v,r;z_0)+\sup\limits_{0<r<R_\star}D^+_2(q,r;z_0)
	\end{equation}
	$$\leq c(M)(D_2^+(q,R_\star;z_0)+E^+(v,R_\star;z_0)+1)\leq$$
	$$\leq  c(M,R_\star,
	\|\nabla v\|_{L_2(Q^+)},\|\nabla q\|_{L_{\frac {12}{11},\frac 32(Q^+)}})=:N.
	$$
Now, having in hands number $M$ and $N$, we may find the number $\varepsilon$ of Proposition \ref{mainlemmabc}.	
	
Let us denote the set of all singulars points of the form $z_0=(x_0,0)$ with $x_0\in b(r_0)$ by $S_b$. Then,  repeating arguments of the proof of Theorem \ref{final} with half balls instead of balls, we  show that the number of elements of $S_b$ is bounded by $M^3/\varepsilon^4$.

Now, it remains to establish that any singular point, belonging to a flat part of the boundary, cannot be the limit point of a sequence of  singular points from the interior of a half ball. To this end,	
	we argue ad absurdum. Let $x_0=x'_0$ with $|x_0|\leq r_0$ be a singular point of $v$ and there exists a sequence  $x^k$ such that $x^k\to x_0$ as $k\to \infty$ and $x^m_3>0$ for all $m$.
	
	Without loss of generality, we may assume $0<x^m_3\leq R_*/2$ for all $m$. In this case, for all $x_*\in B^+(r_0)$ with $0<x_3\leq R_*/2$, the following is valid:
	$$B(x_\star,R_\star)\cap \{x_{3\star}>0 \}\subset B^+(x'_*,2R_\star)\subset B^+.
	$$
Denoting $z_*=(x_*,0)$ and $z'_*=(x'_*,0)$,  observe that 
$$\Theta^+(v,q,r;z'_*):=E^+(v,r;z'_*)+D^+_2(q,r;z'_*)+A^+(v,r;z'_*)+K^+(v,r;z'_*)+$$
$$+C^+(v,r;z'_*)\leq c(M)(E^+(v,2R_*;z'_*)+D^+_2(q,2R_*;z'_*)+$$$$+1)
=:C_1(M,R_*,\|\nabla v\|_{L_2(Q^+)},\|\nabla q\|_{L_{\frac {12}{11},\frac 32}(Q^+)}).$$
for all $0<r\leq R_*$.
 
 It is easy to check that 
$$E(v,x_{3*};z_*)+D_2(v,x_{3*};z_*)\leq c(E^+(v,2x_{*3};z'_*)+D^+_2(q,2x_{*3};z'_*))\leq $$
$$\leq cC_1$$
as $2x_{*3}\leq R_*$. Hence, the number $N$ is determined by  the following inequality
$$E(v,x_{3*};z_*)+D_0(v,x_{3*};z_*)\leq cC_1=:N
$$ 
and  one can find the number $\varepsilon(M,N)$ of Proposition \ref{mainlemma}.

Let us pick up $P$ different elements $x^{k_1}$, $x^{k_2}$,...,$x^{k_P}$ of the sequence $x^k$ assuming that 
$$P>\frac {M^3}{\varepsilon^4}.$$
We let $\gamma=\min\{x^{k_1}_3$, $x^{k_2}_3$,...,$x^{k_P}_3\}>0$
and then select $0<R<\min\{\gamma, R_*/10\}$ so that 
$B(x^{k_i},R)\cap B(x^{k_j},R)=\emptyset$ if $i\neq j$.  Our further arguments are the same  as in the proof of Theorem \ref{final}. Theorem \ref{finalbc} is proved.

\setcounter{equation}{0}
\section{Proof of Proposition \ref{mainlemma}}

We need an auxilary local regularity result, which is in fact  a sufficient condition of  regularity on one scale, see paper \cite{Seregin2016}.

\begin{pro}\label{mainresultinter} Let $v$ and $q$ be a suitable weak solution to the Navier-Stokes equations in $Q(z_0,R)$.

Given $Z>0$, there exist positive numbers $\varepsilon_\star=\varepsilon_*(Z)$ and $c_*=c_*(Z)$ such that if two conditions
$$\frac 1{R^2}\int\limits_{Q(z_0,R)}|v|^3dxdt< \varepsilon_*(Z)$$ 
and
$$D_0(q,R;z_0)=\frac 1{R^2}\int\limits_{Q(z_0,R)}|q-[q]_{B(x_0,R)}|^\frac 32dxdt<Z$$
hold, then  $v$ and $\nabla v$ are H\"older continuous is the closure of  $Q(z_0,R/2)$. Moreover,
$$\sup\limits_{z\in Q(z_0,R/2)}|v(z)|+|\nabla v(z)|\leq \frac {c_*(Z)}R.$$
\end{pro}


As in paper \cite{CWY}, we argue as absurdum. Indeed, if we assume that the statement of Proposition \ref{mainlemma} is false, then, according to Proposition \ref{mainresultinter}, there are positive numbers $M$ and $N$ such that there exist a sequence of suitable weak solutions $v^k$ and $q^k$, sequences of numbers $0<r_k\leq 1/2$ and $\varepsilon_k\to +0$ with the following properties:
\begin{equation}
	\label{wl3}
	\sup\limits_k\|v^k\|_{L_\infty(-1,0;L^{3,\infty}(B))}\leq M;
\end{equation} 
\begin{equation}
	\label{startingpoint}
\sup\limits_kD_0(q^k,1)\leq \sup\limits_k(D_0(q^k,1)+E(v^k,1))\leq N;	
\end{equation}
\begin{equation}
	\label{maincondk}
	\frac 1{r_k^3}|\{x\in B(r_k):\,\,|v^k(x,0)|>\frac {\varepsilon_k}{r_k}\}|\leq\varepsilon_k
\end{equation}
for all $k=1,2,...$;
\begin{equation}
	\label{sing}
	\frac 1{\varrho^2}\int\limits_{Q(\varrho)}|v^k|^3dz>\frac 12 \varepsilon_*(N)
\end{equation}
for all $\varrho\in [2\varepsilon_kr_k,r_k]$.

Moreover, the same arguments as in the proof of the main theorem  lead to the inequality for energy scale invariant quantities:
\begin{equation}
	\label{scale inv}
\Theta(v^k,q^k,r;z_0):=	A(v^k,r;z_0)+C(v^k,r;z_0)+E(v^k,r;z_0)+
	\end{equation}
$$+K(v^k,r;z_0)+D_0(q^k,r;z_0)\leq c(M)(N+1)$$
for all $0<r\leq 1/2$ and for all $z_0\in Q(1/2)$.

Now,  our functions can be scaled  in the following way:
$$u^k(y,s)=r_kv^k(x,t),\qquad p^k(y,s)=r^2_kq^k(x,t),$$
where $x=r_ky$, $t=r_k^2s$ and $e=(y,s)\in Q(1/r_k)$. New functions $u^k$ and $p^k$ satisfy the Navier-Stokes equations in  $Q(1/r_k)$,  
$$
\sup\limits_k(D_0(p^k,1/r_k)
+E(u^k,1/r_k))<N,$$
and 
$$\|u^k\|_{L_\infty(-1/r^2_k,0;L^{3,
\infty}(B(1/r_k))}\leq M.
$$

Without loss of generality, one may assume  that $r_k\to r_*$ as $k\to \infty$. There are two case: $r_*=0$ and $r_*>0$.

Let us first consider the case $r_*=0$. Here,  we can  fix an arbitrary space-time point $e_0=(y_0,s_0)$, a number $0<R<\frac 12\frac 1{r_k}$ and make change of variables in (\ref{scale inv}) in order to get
\begin{equation}
	\label{scaled scale inv}
	\Theta(u^k,p^k,R;e_0)\leq c(M)(N+1).
\end{equation}
for sufficiently large $k$. Moreover, (\ref{maincondk}) and  (\ref{sing}) can be transformed into the following:
\begin{equation}
	\label{scalemaincond}
	|\{y\in B:\,|u^k(y,0)|>\varepsilon_k\}|\leq\varepsilon_k
\end{equation}
and 
\begin{equation}
	\label{scalesing}
\frac 1{\varrho^2}\int\limits_{Q(\varrho)}|u^k|^3de>\frac 12 \varepsilon_*(N)
\end{equation}
for all $\varrho\in [2\varepsilon_k,1]$ and for all $k$.	

Higher derivatives can be evaluated as in \cite{Seregin}. So,
\begin{equation}
	\label{Higherder}
	\Sigma(u^k,p^k,R):=
	\frac 1R\Big[\|\partial_tu^k\|_{L_{\frac 98,\frac 32}(Q(R))}+\|\nabla^2 u^k\|_{L_{\frac 98,\frac 32}(Q(R))}+\end{equation}
$$
+\|\nabla p^k\|_{L_{\frac 98,\frac 32}(Q(R))}\Big]\leq c\Big[A^\frac 13(u^k,2R)E^\frac 23(u^k,2R)+A^\frac 12(u^k,2R)+$$
$$+E^\frac 12(u^k,2R)+D_0^\frac 23(p^k,2R)\Big]$$
for all $0<R<\frac 14\frac 1{r_k}$.

Now, let us pass to the limit as $k\to\infty$, taking into account estimates (\ref{scaled scale inv}) and (\ref{Higherder}). Then,
after using known compactness  arguments, we get the so-called  local energy ancient solution $u$ and $p$, having the following properties:
$$u^k\to u$$ 
in $L_3(Q(R))$ and in $C([-R^2,0];L_\frac 98(B(R)))$ for any $R>0$;
$$p^k\rightharpoonup p$$
in $L_\frac 32 (Q(R))$ for any $R>0$;

the pair $u$ and $p$ is a suitable weak solutionto the Navier-Stokes equations in each $Q(R)$; 
$$\|u\|_{L_\infty(-\infty,0;L^{3,\infty}(\mathbb R^3))}\leq M;$$
$$\Theta(u,p,R;e_0)\leq c(M)(N+1)$$
for any $R>0$ and $e_0\in Q_-$;
$$\Sigma(u,p,R)\leq c(M,N)$$
for any $R>0$;
$$\frac 1{\varrho^2}\int\limits_{Q(\varrho)}|u|^3de\geq \frac 12 \varepsilon_*(N)$$
for any $\varrho\in ]0,1]$, and finally
\begin{equation}
\label{zeroinball}	u(x,0)=0
\end{equation}
for any $x\in B$.  

The latter identity follows from the known inequality
$$|\{y\in B: |u(y,0)|>\alpha\}|\leq |\{y\in B: |u^k(y,0)|>\alpha/2\}|+$$$$+|\{y\in B: |u^k(y,0)-u(y,0)|>\alpha/2\}|$$
that is valid for any $\alpha>0$.

Now, let us consider the case $r_*>0$. Our first remark is that (\ref{scaled scale inv}) remains to be true for all $e_0=(y_0,s_0)$  from the unit parabolic ball  $ Q$ and for the same $R$. Moreover, relationships (\ref{scalemaincond})-(\ref{Higherder}) are completely the same as well. Repeating the same compactness arguments, we can easily pass to the limit as $k\to\infty$ and conclude that there exist functions $u$ and $p$ with the following properties:

the pair $u$ and $p$ is a suitable weak solutionto the Navier-Stokes equations in each $Q(1/4)$; 
$$\|u\|_{L_\infty(-1/4^2,0;L^{3,\infty}(B(1/4))}\leq M;$$
$$E(u,1/4)+D_0(p,1/4)\leq c(M)(N+1);$$
$$\frac 1{\varrho^2}\int\limits_{Q(\varrho)}|u|^3de\geq \frac 12 \varepsilon_*(N)$$
for any $\varrho\in ]0,1/4]$, and finally
$u(y,0)=0
$
for any $y\in B(1/4)$. 

Obviously, the restriction of $u$ and $p$ of the first case $r_*=0$ to the parabolic ball $Q(1/4)$ have the properties as above and in what follows we shall work with such a restriction.

The crucial point here is a reduction to backward uniqueness for the heat operator with lower order terms, see \cite{ESS2003}.
To this end, we 
select  a sequence of positive numbers $\varrho_k$, tending to zero. Then, the new scaling is:
$$U^k(y,s)=\varrho_ku(x,t),\qquad P^k(y,s)=\varrho^2_kp(x,t)$$
where $x=\varrho_ky$, $t=\varrho^2_ks$.  Repeating arguments of the first part of the proof, we find the following relationships:

given $e_0\in Q_-$,
$$\Theta(U^k,P^k,R;e_0)+\Sigma(U^k,P^k,R/2)\leq c(M,N)$$
for any $0<R<1/(2\varrho_k)$  and for sufficiently large $k$;
$$\|U^k\|_{L_\infty(-1/\varrho_k^2,0;L^{3,\infty}(B(1/\varrho_k))}\leq M;$$

$$U^k(x,0)=0
$$
for all $x\in B(1/(4\varrho_k))$;
$$\frac 1{(\varrho/\varrho_k)^2}\int\limits_{Q(\varrho/\varrho_k)}|U^k(z)|^3dz\geq \frac 12\varepsilon_*(N)$$
for all $0<\varrho\leq1/4$. 

Let $\varrho=\varrho_k$ and $k$ tend to infinity and let us see what happens.
The same arguments as in the first scaling lead to the following: there exists a local energy ancient solution $w$ with the associated pressure $r$ such that:

$$\Theta(w,r,\varrho;z_0)+\Sigma(w,r,\varrho)\leq c(M,N)$$
for any $\varrho>0$ and for any $z_0\in Q_-$;
$$\|w\|_{L_\infty(-\infty,0;L^{3,\infty}(\mathbb R^3))}\leq M;$$
$$w(x,0)=0
$$
for all $x\in\mathbb R^3$;
$$
\int\limits_{Q
}|w|^3dz\geq \frac 12\varepsilon_*(N)>0.$$
In order to apply the approach based on the backward uniqueness, we need to show that solution $w$ has a certain decay at infinity, for example, to prove that  $w$ and $\nabla w$ belong to $L_\infty((\mathbb R^3\setminus B(R))\times ]-2,0[)$ for some $R>0$.
Just for completeness, we repeat arguments from the paper \cite{CWY}. Indeed, by the definition of weak Lebesgue spaces, we find that 
$$|\{(x,t)\in \mathbb R^3\times ]-3,0[:\,\,|w(x,t)|>\gamma\}|\leq \frac 3{\gamma^3}M^3<\infty.$$
Hence, for any positive number $\eta>0$, there exists $R=R(\gamma)>0$ such that 
$$|\{(x,t)\in(\mathbb R^3\setminus B(R(\gamma)))\times ]-3,0[:\,\,|w(x,t)|>\gamma\}|<\eta.
$$
So, if $Q(z_0,1)\in (\mathbb R^3\setminus B(R(\gamma)))\times ]-3,0[$, then we have 
$$D_0(r,1;z_0)
\leq c(M,N)=Z $$
and 
$$C(w,1;z_0)\leq \gamma^3|Q(z_0,1)|+\int\limits_{\{(x,t)\in Q(z_0,1):\,|w(x,t)|>\gamma\}}|w|^3dz\leq 
$$
$$\leq c\gamma^3+\Big(\int\limits_{Q(z_0,1)}|w|^4dz\Big)^\frac 34|\{(x,t)\in Q(z_0,1):\,|w(x,t)|>\gamma\}|^\frac 14\leq $$
$$\leq c\gamma^3+K^\frac 34(w,1;z_0)\eta^\frac 14 \leq c\gamma^3+c(M,N)\eta^\frac 14.$$
We select first $\gamma$ and then $\eta$ so that the right hand side of the latter inequality is   less that $\varepsilon_\star(Z)$. Then, one can conclude, see Proposition \ref{mainresultinter}, that, for any $z_0\in (\mathbb R^3\setminus B(\eta))\times ]-1,0[$, 
$$|u(z_0)|+|\nabla u(z_0)|\leq c_\star(Z). $$

Now, using arguments of the paper \cite{ESS2003}, we show that $w\equiv 0$ in $\mathbb R^3\times ]-1,0[$. This is a contradiction.  So, Proposition \ref{mainlemma} is proved.

\setcounter{equation}{0}
\section{Proof of Proposition \ref{mainlemmabc}}

Since our proof of the proposition is  similar to the proof of Proposition \ref{mainlemma}, we just outline it. We start with 
 a certain boundary regularity condition, following   the paper \cite{Seregin2016}.

\begin{pro}\label{mainresultbc} Let $v$ and $q$ be a suitable weak solution to the Navier-Stokes equations in $Q^+(z_0,R)$ in the sense of Definition \ref{swsbc}.
Given $Z>0$, there exist positive numbers $\varepsilon_\star=\varepsilon_*(Z)$ and $c_*=c_*(Z)$ such that if two conditions
$$\frac 1{R^2}\int\limits_{Q^+(z_0,R)}|v|^3dxdt< \varepsilon_*(Z)$$ 
and
$$\frac 1{R^2}\int\limits_{Q^+(z_0,R)}|q-[q]_{B^+(x_0,R)}|^\frac 32dxdt<Z$$
hold, then  $v$ is H\"older continuous is the closure of  $Q^+(z_0,R/2)$. Moreover,
$$\sup\limits_{z\in Q^+(z_0,R/2)}|v(z)|\leq \frac {c_*(Z)}R.$$
\end{pro}

Assume that Proposition \ref{mainlemmabc} is false. Then, there exist sequeneces $v^k$, $q^k$, $0<r_k\leq 1/2$, and $\varepsilon_k\to+0$ such that
$$ \sup\limits_kD_0^+(q^k,1)\leq c \sup\limits_kD_2^+(q^k,1)\leq cN,
$$
$$
\sup\limits_k\|v^k\|_{L_\infty(-1,0;L^{3,\infty}(B^+))}\leq M,
$$
$$
\frac 1{r_k^3}|\{x\in B^+(r_k):\,|v^k(x,0)|>\frac {\varepsilon_k}{r_k}\}|<\varepsilon_k,
$$
but
$$\frac 1{\varrho^2}\int\limits_{Q^+(\varrho)}|v^k|^3dz>\frac 12\varepsilon_*(cN)
$$
for all $\varrho\in [2\varrho_kr_k,r_k]$.
 
 Then, we have a typical estimate of certain energy scale invariant quantities:
 $$A^+(v^k,r;z_0)+C^+(v^k,r;z_0)+E^+(v^k,r;z_0)+K^+(v^k,r;z_0)+ $$
$$+D^+2(q^k,r;z_0)\leq C(M)(N+1)$$
for all $z_0\in Q^+(1/2)$ zuch that $z_0=(x_0,t_0)$ and $x_0=x_0'$.

We let $\omega(x_0,r)=B(x_0,r)\cap \mathbb R^3_+$ and $Q_\omega(z_0,r)=\omega(x_0,r)\times ]t_0-r^2,t_0[$.
 Using scaling arguments, we get the main estimate
 $$\Theta_\omega(v^k,q^k,r;z_0):=A_\omega(v^k,r;z_0)+C_\omega(v^k,r;z_0)+E_\omega(v^k,r;z_0)+$$$$+K_\omega(v^k,r;z_0)+D_{2\omega}(q^k,r;z_0) \leq c(M)(N+1)$$
 for all $0<r\leq \frac 14$ and $z_0\in Q^+(1/4)$.
 
Next, we do scaling $u^k(y,s)=r_kv^k(x,s)$,  $p^k(y,s)=r^2_kq^k(x,s)$, where $x=r_ky$, $t=r^2_ks$ and  $e=(y,s)\in Q^+(1/r_k)$. Here, we are going to consider the only case in which $r_k\to0$ as $k\to \infty$, leaving the second case to the reader.

 From the previous estimates, one can  deduce the following:
$$\sup\limits_k(E^+(u^k,1/r_k)+D_2^+(p^k,1/r_k))\leq N.$$
 Moreover, we fix $e_0=(y_0,s_0)$ and, for $0<R<1/(4r_k)$, find
 $$\Theta_\omega(u^k,p^k,R;e_0)\leq c(M)(N+1),$$ 
  $$|\{y\in B^+:\,\,|u^k(y,0)|>\varepsilon_k\}|<\varepsilon_k,$$
and 
 $$\frac 1{\varrho^2}\int\limits_{Q^+(\varrho)}|u^k|^3de>\frac 12\varepsilon_*(cN).$$
In order to provide compactness, higher derivatives are evaluated so that:
$$\Sigma^+(u^k,p^k,R)=\frac 1{R^{\frac {13}{12}}}[\|\partial_tu^k\|_{L_{\frac {12}{11},\frac 32,Q^+(R)}}+\|\nabla^2 u^k\|_{L_{\frac {12}{11},\frac 32,Q^+(R)}}+$$$$+\|\nabla p^k\|_{L_{\frac {12}{11},\frac 32,Q^+(R)}}]\leq $$
$$\leq c[D_0^+(p^k,2R)+ (C^+)^\frac 13(u^k,2R)+(E^+)^\frac 12(u^k,2R)+$$ 
 $$+(A^+)^\frac 14(u^k,2R)(E^+)^\frac 12(u^k,2R) (C^+)^\frac 16(u^k,2R)]
 $$
 for $0<R<1/(8r_k)$.
 
Passing to the limit as $k\to\infty$, we have (without loss of generality) the following:
$u^k\to u$  in $L_3(Q^+(R))$ and in $C([-R^2,0];L_\frac {12}{11}(B^+(R)))$ and 
$p^k\rightharpoonup p$ in $L_\frac 32(Q^+(R))$ for all $R>0$. Let us list the properies of the limit pair $u$ and $p$. It is a suitable weak solution to the Navier-Stokes equations in each $Q^+(R)$;
$$\|u\|_{L_\infty(-\infty,0;L^{3,\infty}(\mathbb R^3_+))}\leq M;$$
$$\Theta_\omega(u,p,R;e_0)\leq c(M)(N+1)$$
 for all $R>0$ and for all $e_0\in Q^+_-:=\mathbb R^3_+\times ]-\infty,0[$;
 $$\Sigma^+(u,p,R)\leq c(M,N)$$
 for all positive $R$;
 $$\frac 1{\varrho^2}\int\limits_{Q^+(\varrho)}|u|^3dz\geq\frac 12\varepsilon_*(cN)$$
 for $0<\varrho\leq 1$;
 $$u(x,0)=0$$
 for $x\in B^+(0)$.

 Our further arguments are the same in the proof of Proposition \ref{mainlemma}. We need to adopt the  last part of the proof based on the backward uniqueness since the boundary local regularity, in general, does not provide boundedness of $\nabla u$ up to the boundary. To this end, let us fix a positive number $h$ and apply the interior boundary regularity result (as we did in the previous section) in order to show that $\nabla u$ is bounded in $(R^3_++he_3)\setminus B(R_0)$ for a large number $R_0$. Then we can use  known arguments, based on the backward uniqueness and the unique continuation through spatial boundaries. This implies $u\equiv 0$ in $R^3_++he_3$ for all $h>0$ and thus we arrive at the contardiction. The proposition is proved.

\end{document}